\documentclass{gtart_h}

\def\ifplaintex{\expandafter\ifx\csname documentclass\endcsname\relax}

\def\gtp{{\mathsurround=0pt\it $\cal G\mskip-2mu$eometry \&\ 
$\cal T\!\!$opology $\cal P\!$ublications}}  

\def\recd{{\small Received:\qua\receiveddate\ifx\reviseddate\relax
\else\qquad Revised:\qua\reviseddate\fi\par}} 

\def\lognumber#1{\def\thelognumber{#1}}
\def\volumenumber#1{\def\thevolumenumber{#1}}
\def\volumeyear#1{\def\thevolumeyear{#1}}
\def\papernumber#1{\def\thepapernumber{#1}}
\def\pagenumbers#1#2{\def\startpage{#1}\def\finishpage{#2}}
\def\published#1{\def\publishdate{#1}}

\def\received#1{\def\receiveddate{#1}}
\def\revised#1{\def\reviseddate{#1}}
\def\accepted#1{\def\accepteddate{#1}}

\def\asciiaddress#1{\def\theasciiaddress{#1}}
\def\asciiemail#1{\def\theasciiemail{#1}}

\long\def\asciiabstract#1{\long\def\theasciiabstract{#1}}
\def\asciikeywords#1{\def\theasciikeywords{#1}}

\let\\\par\let\thelognumber\relax\let\thevolumenumber\relax
\let\thepapernumber\relax\let\thevolumeyear\relax\let\startpage\relax
\let\finishpage\relax\let\publishdate\relax\let\receiveddate\relax
\let\reviseddate\relax\let\accepteddate\relax\let\theasciititle\relax
\let\theasciiauthors\relax\let\theasciiaddress\relax
\let\theasciiabstract\relax\let\theasciikeywords\relax

\let\theasciiemail\relax

\ifplaintex
\font\logobig=cmssbx10 scaled 3836
\font\logomed=cmssbx10 scaled 2557
\else
\font\logobig=cmssbx10 scaled 4200
\font\logomed=cmssbx10 scaled 2800
\fi

\long\def\makeagttitle{   
\count0=\startpage
\agt\hfill      
\hbox to 45truept{\vbox to 0pt{\vglue -13truept{\logomed A\kern -.37em{\logobig 
T}\kern -.38em G}\vss}\hss}
\break
{\small Volume \thevolumenumber\ (\thevolumeyear)
\startpage--\finishpage\nl
Published: \publishdate}

\vglue .25truein

{\parskip=0pt\leftskip 0pt plus
1fil\def\\{\par\smallskip}{\Large\bf\thetitle}\par\medskip} \vglue
0.05truein

{\parskip=0pt\leftskip 0pt plus 1fil\def\\{\par}{\sc\theauthors}
\par\medskip}%
 
\vglue 0.03truein 

{\small\leftskip 25truept\rightskip 25truept{\bf Abstract}\stdspace\theabstract

{\bf AMS Classification}\stdspace\theprimaryclass
\ifx\thesecondaryclass\relax\else; \thesecondaryclass\fi\par
{\bf Keywords}\stdspace \thekeywords\par}\vglue 7truept

}   

\ifplaintex
\font\phead=cmsl9 scaled 950
\font\pnum=cmbx10 scaled 913
\font\pfoot=cmsl9 scaled 950
\headline{\vbox to 0pt{\vskip -4.5mm\line{\small\phead\ifnum
\count0=\startpage ISSN 1472-2739 (on-line) 1472-2747 (printed)
\hfill {\pnum\folio}\else\ifodd\count0\def\\{ }%
\ifx\theshorttitle\relax\thetitle\else\theshorttitle\fi\hfill{\pnum\folio}
\else\def\\{ and }{\pnum\folio}\hfill\ifx\theshortauthors\relax\theauthors
\else\theshortauthors\fi\fi\fi}\vss}}
\footline{\vbox to 0pt{\vglue 0mm\line{\small\pfoot\ifnum\count0=\startpage
\copyright\ \gtp\hfill\else
\agt, Volume \thevolumenumber\ (\thevolumeyear)\hfill\fi}\vss}}
\else
\font=cmsl9 scaled 1050
\font\lnum=cmbx10 
\font=cmsl9 scaled 1050
\makeatletter
\def\@oddhead{{\small\lhead\ifnum\count0=\startpage ISSN 1472-2739 
(on-line) 1472-2747 (printed)\hfill {\lnum\number\count0}\else\ifodd\count0
\def\\{ }\ifx\theshorttitle\relax \thetitle \else\theshorttitle\fi\hfill
{\lnum\number\count0}\else\def\\{ and }{\lnum\number\count0}
\hfill\ifx\theshortauthors\relax 
\theauthors\else\theshortauthors\fi\fi\fi}}\def\@evenhead{\@oddhead}
\def\@oddfoot{\small\lfoot\ifnum\count0=\startpage\copyright\ \gtp\hfill\else
\agt, Volume \thevolumenumber\ (\thevolumeyear)\hfill\fi}
\def\@evenfoot{\@oddfoot}
\makeatother
\fi
\let\maketitlepage\makeagttitle
\let\makeshorttitle\maketitlepage
\let\maketitle\maketitlepage

\newwrite\gtoutfile
\long\gdef\makeheadfile{  
{\def\\{, }\def\s{ }
\immediate\openout\gtoutfile head.xxx
\immediate\write\gtoutfile{Proxy-for: \ifx\theasciiauthors\relax
\theauthors\else\theasciiauthors\fi\s<\ifx\theasciiemail\relax\theemail\else\theasciiemail\fi>}
\immediate\write\gtoutfile{\noexpand\\}
\immediate\write\gtoutfile{Authors: \ifx\theasciiauthors\relax
\theauthors\else\theasciiauthors\fi}
{\def\\{ }\immediate\write\gtoutfile{Title: \ifx\theasciititle\relax
\thetitle\else\theasciititle\fi}}
\immediate\write\gtoutfile{Subj-class: GT or SG, GR etc}
\immediate\write\gtoutfile{MSC-class: \theprimaryclass\ifx\thesecondaryclass\relax\else, \thesecondaryclass\fi}
\immediate\write\gtoutfile{Journal-ref: Algebraic and Geometric Topology \thevolumenumber\s
(\thevolumeyear) \startpage-\finishpage}
\immediate\write\gtoutfile{Comments: Published by Algebraic and
Geometric Topology at}
\immediate\write\gtoutfile{\s\s\s  http://www.maths.warwick.ac.uk/agt/AGTVol\thevolumenumber/agt-\thevolumenumber-\thepapernumber.abs.html}
\immediate\write\gtoutfile{\noexpand\\}
\immediate\write\gtoutfile{}
\ifx\theasciiabstract\relax
\immediate\write\gtoutfile{\theabstract}\else
\immediate\write\gtoutfile{\theasciiabstract}\fi
\immediate\write\gtoutfile{}
\immediate\write\gtoutfile{\noexpand\\}
\immediate\write\gtoutfile{}
\immediate\closeout\gtoutfile}}  

\def\maketitlepage{\makeagttitle\makeheadfile}
\let\makeshorttitle\maketitlepage
\let\maketitle\maketitlepage

\lognumber{4}
\volumenumber{4}
\volumeyear{2004}
\papernumber{4}
\published{7 February 2004}
\pagenumbers{49}{71}
\received{21 November 2002}
\revised{15 October 2003}
\accepted{26 November 2003}

\usepackage{graphicx}
\nocolon

\newtheorem{theorem}{Theorem}
\newtheorem{lemma}{Lemma}

\begin{document}

\title{The boundary-Wecken classification of surfaces}

\author{Robert F. Brown\\Michael R. Kelly} 
\address{Department of Mathematics, University of California\\Los Angeles, CA 90095-1555, USA} 
\gtemail{\mailto{rfb@math.ucla.edu}\qua{\rm and}\qua\mailto{kelly@loyno.edu}}

\secondaddress{Department of Mathematics and Computer Science, Loyola University\\6363 St. Charles Avenue, New Orleans, LA 70118, USA} 

\asciiaddress{Department of Mathematics, University of California\\Los Angeles, CA 90095-1555, USA\\and\\Department of Mathematics and Computer Science, Loyola University\\6363 St. Charles Avenue, New Orleans, LA 70118, USA}
\asciiemail {rfb@math.ucla.edu, kelly@loyno.edu}

\begin{abstract}
Let $X$ be a compact $2$-manifold with nonempty boundary $\partial X$
and let $f \co (X, \partial X) \to (X, \partial X)$ be a
boundary-preserving map.  Denote by $MF_{\partial}[f]$ the minimum
number of fixed point among all boundary-preserving maps that are
homotopic through boundary-preserving maps to $f$.  The relative
Nielsen number $N_{\partial}(f)$ is the sum of the number of essential
fixed point classes of the restriction $\bar f \co \partial X \to
\partial X$ and the number of essential fixed point classes of $f$
that do not contain essential fixed point classes of $\bar f$.  We
prove that if $X$ is the M\"obius band with one (open) disc removed,
then $MF_{\partial}[f] - N_{\partial}(f) \le 1$ for all maps $f \co
(X, \partial X) \to (X, \partial X)$.  This result is the final step
in the boundary-Wecken classification of surfaces, which is as
follows.  If $X$ is the disc, annulus or M\"obius band, then $X$ is
boundary-Wecken, that is, $MF_{\partial}[f] = N_{\partial}(f)$ for all
boundary-preserving maps.  If $X$ is the disc with two discs removed
or the M\"obius band with one disc removed, then $X$ is not
boundary-Wecken, but $MF_{\partial}[f] - N_{\partial}(f) \le 1$.  All
other surfaces are totally non-boundary-Wecken, that is, given an
integer $k \ge 1$, there is a map $f_k \co (X, \partial X) \to (X,
\partial X)$ such that $MF_{\partial}[f_k] - N_{\partial}(f_k) \ge k$.
\end{abstract}

\asciiabstract{%
Let X be a compact 2-manifold with nonempty boundary dX and let f: (X,
dX) --> (X, dX) be a boundary-preserving map.  Denote by MF_d[f] the
minimum number of fixed point among all boundary-preserving maps that
are homotopic through boundary-preserving maps to f.  The relative
Nielsen number N_d(f) is the sum of the number of essential fixed
point classes of the restriction f-bar : dX --> dX and the number of
essential fixed point classes of f that do not contain essential fixed
point classes of f-bar.  We prove that if X is the Moebius band with
one (open) disc removed, then MF_d[f] - N_d(f) < 2 for all maps f :
(X, dX) --> (X, dX).  This result is the final step in the
boundary-Wecken classification of surfaces, which is as follows.  If X
is the disc, annulus or Moebius band, then X is boundary-Wecken, that
is, MF_d[f] = N_d(f) for all boundary-preserving maps.  If X is the
disc with two discs removed or the Moebius band with one disc removed,
then X is not boundary-Wecken, but MF_d[f] - N_d(f) < 2.  All other
surfaces are totally non-boundary-Wecken, that is, given an integer k
> 0, there is a map $f_k : (X, dX) --> (X, dX) such that MF_d[f_k] -
N_d(f_k) >= k.}



\asciikeywords{Boundary-Wecken, relative Nielsen number, punctured
Moebius band, boundary-preserving map}

\maketitle

\setcounter{section}{-1}

\section{Introduction}

Given a map $f \co X \to X$ of a compact $n$-manifold with
(possibly empty) boundary to itself, let $N(f)$ denote the Nielsen
number of $f$, that is, the number of essential fixed point classes.
For $MF[f]$ the minimum number of fixed points among all
maps homotopic to $f$, we have $N(f) \le MF[f]$.  Wecken proved in
\cite{w} that if $X$ is an $n$-manifold with $n \ge 3$, then $N(f) =
MF[f]$ for all maps $f \co X \to X$.  Consequently, the property:
$N(f) = MF[f]$ for every map $f \co X \to X$ has come to be called
the {\it Wecken property} of a manifold $X$ and a manifold with this
property is said to be a {\it Wecken manifold}.  If $X$ is a
$2$-manifold with Euler characteristic non-negative, then
$X$ is a Wecken manifold.  The Wecken property of the Klein bottle is
a consequence of Theorem 5.1 of \cite{h} (see Corollary 8.3 of 
\cite{hkw}),
the property is verified for the torus in \cite{b} and
for the projective plane in \cite{j1}.  The verification of the
Wecken property for the four remaining such $2$-manifolds (sphere,
disc, annulus, M\"obius band) is easy.

On the other hand, if $X$ is a $2$-manifold with negative Euler
characteristic, then it is proved in \cite{j2} that $X$ is {\it totally
non-Wecken}, that is, given an integer $k \ge 1$, there is a map $f_k
\co X \to X$ such that $MF[f_k] - N(f_k) \ge k$.   Thus we have

\begin{theorem}[Wecken classification of surfaces] Let $X$
be a compact $2$-man\-ifold.  

{\rm(a)}\qua If the Euler characteristic of $X$
is non-negative, then $X$ is
Wecken.  

{\rm(b)}\qua All other surfaces are totally non-Wecken.
\end{theorem}

Now let $X$ denote a compact $n$-manifold with nonempty boundary
$\partial X$ and consider boundary-preserving self-maps $f \co (X,
\partial X) \to (X, \partial X)$.  We write the restriction of $f$ to
the boundary as $\bar f \co \partial X \to \partial X$.  The
minimum number of fixed points among all maps of pairs homotopic to
$f$ through boundary-preserving maps is written as $MF_\partial[f]$.  
The
relative Nielsen number, that we denote by $N_\partial(f)$, is defined
to be the sum of the number of essential fixed point classes of $\bar
f$ and the number of essential fixed point classes of $f$ that do not
contain essential fixed point classes of $\bar f$.  We always have
$N_\partial(f) \le MF_\partial[f]$ and the $n$-manifold $X$ is said
to be {\it boundary-Wecken} if $N_\partial(f) = MF_\partial[f]$ for
all maps of pairs $f \co (X, \partial X) \to (X, \partial X)$.
Schirmer \cite{s} proved that all $n$-manifolds with nonempty
boundary are boundary-Wecken if $n \ge 4$.

A compact $n$-manifold $X$ with nonempty boundary is said to be {\it 
totally
non-boundary-Wecken} if, given an integer $k \ge 1$, there is a map
$f_k \co (X, \partial X) \to (X, \partial X)$ such that
$MF_\partial[f_k] - N_\partial(f_k) \ge k$.  The
manifold $X$ is said to be {\it almost boundary-Wecken} if $X$ is
not boundary-Wecken but there exists an integer $B \ge 1$ such that
$MF_\partial[f] - N_\partial(f) \le B$ for all maps $f \co (X,
\partial X) \to (X, \partial X)$.

The purpose of this paper is to prove the following theorem:

\begin{theorem}
Let $X$ be the M\"obius band
with one (open) disc removed and let
$f \co (X, \partial X) \to (X, \partial X)$ be any
map, then $MF_\partial[f] - N_\partial(f) \le 1$.
\end{theorem}

This result will allow us to complete the proof of the following:

\begin{theorem}[Boundary-Wecken classification of surfaces]
Let $X$ be a compact $2$-manifold with nonempty boundary.

{\rm(a)}\qua If $X$ is the disc, annulus or M\"obius band,
then $X$ is boundary-Wecken.  

{\rm(b)}\qua If $X$ is the disc with two discs removed or the M\"obius
band with one disc removed, then $X$ is almost boundary-Wecken, with
$B = 1$.

{\rm(c)}\qua All other surfaces with boundary are totally
non-boundary-Wecken.
\end{theorem}

\proof Conclusion (a) is proved in \cite{bs}.  With regard to
conclusion (b), the proof that the disc with two discs removed
is not boundary-Wecken is in \cite{k1} and it is proved in \cite{bs}
that it is almost boundary-Wecken with $B = 1$.  It is shown in 
\cite{k2}
that the M\"obius band with one disc removed is not boundary-Wecken so 
Theorem
2 will complete the proof of conclusion (b).   The proof that the
disc with three or more discs removed is totally
non-boundary-Wecken is in \cite{n}.  For the M\"obius band with
two or more discs removed, the proof that it is totally
non-boundary-Wecken is in \cite{k2}.  For $X$ any surface not yet 
listed,
results from \cite{bs} and \cite{k1} imply that $X$ is
totally non-boundary-Wecken.  Thus we have conclusion (c).
\endproof

To begin the proof of Theorem 2, we show in Section 1 that, for $X$
the M\"obius band with one disc removed, any map
$f \co (X, \partial X) \to (X, \partial X)$ can be homotoped as a
map of pairs
to a convenient ``reduced form".  For a map in that form, the inverse
images of the arcs where handles are attached to a disc to form $X$
are properly embedded $1$-dimensional submanifolds of $X$.  We
classify the possible such submanifolds in Section 1 also.

Now let $C_1$ and $C_2$ be the components of
$\partial X$ and let $f_i \co C_i \to \partial X$ be
the restriction of $f$ to $C_i$ for $i = 1, 2$.   The following
terminology will be convenient for presenting the various cases to be
checked in order to prove Theorem 2.  If both the $f_i$ are essential
maps, we will say that $f$ is {\it boundary
essential}.
If both $f_i$ are inessential maps, then $f$ is {\it boundary
inessential}.
If one of the maps $f_1$ or $f_2$ is essential and the other
is inessential, we say that
$f$ is {\it boundary semi-essential}.  Also, as in \cite{bs}, we let
$Im_\partial(f)$ denote the number of components of $\partial X$ that
contain points of $f(\partial X)$.

In Section 2, we
prove Theorem 2 in the case that $f$ is boundary inessential.
Section 3 is devoted to the proof of Theorem 2
when $f$ is boundary essential.  We complete the proof of Theorem 2 in 
Section
4, by proving it for boundary semi-essential maps.

Brian Sanderson helped us with the proof of Therem 8.  We also thank 
the referee who gave us many constructive suggestions that have
significantly improved the exposition.

\section{Reduced form}

{\bf Notation}\qua In Figure 1, we represent the punctured M\"obius band
$X$ as a union of discs $X = D \cup H_1 \cup H_2$ where the handle
$H_1$ is twisted,
that is, $D \cup H_1$ is a M\"obius band, and the handle $H_2$ is not
twisted: $D \cup H_2$ is an annulus.  We have $D \cap H_1 = A_1 \cup
A_2$ and $D \cap H_2 = A_3 \cup A_4$ and we let $A$ be the union of
all the arcs $A_j$.   We denote the arc that is the intersection of
$\partial D$ with the component $C_2$ of $\partial X$ by $T_2$ and
the three arcs that are the components of $\partial D\cap C_1$ are
called $T_1, T_3$ and $T_4$, as the figure shows.  We choose points
$x_1 \in T_1$ and $x_2 \in T_2$ as well as $x_0 \in int\, X$ on the
line segment, that we denote by $[x_1, x_2]$, with endpoints $x_1$
and $x_2$.  The simple closed curves $h_1$ and $h_2$ generate the
fundamental group of $X$ based at $x_0$.

\begin{figure}[ht!]
\centering
\includegraphics[scale=0.4]{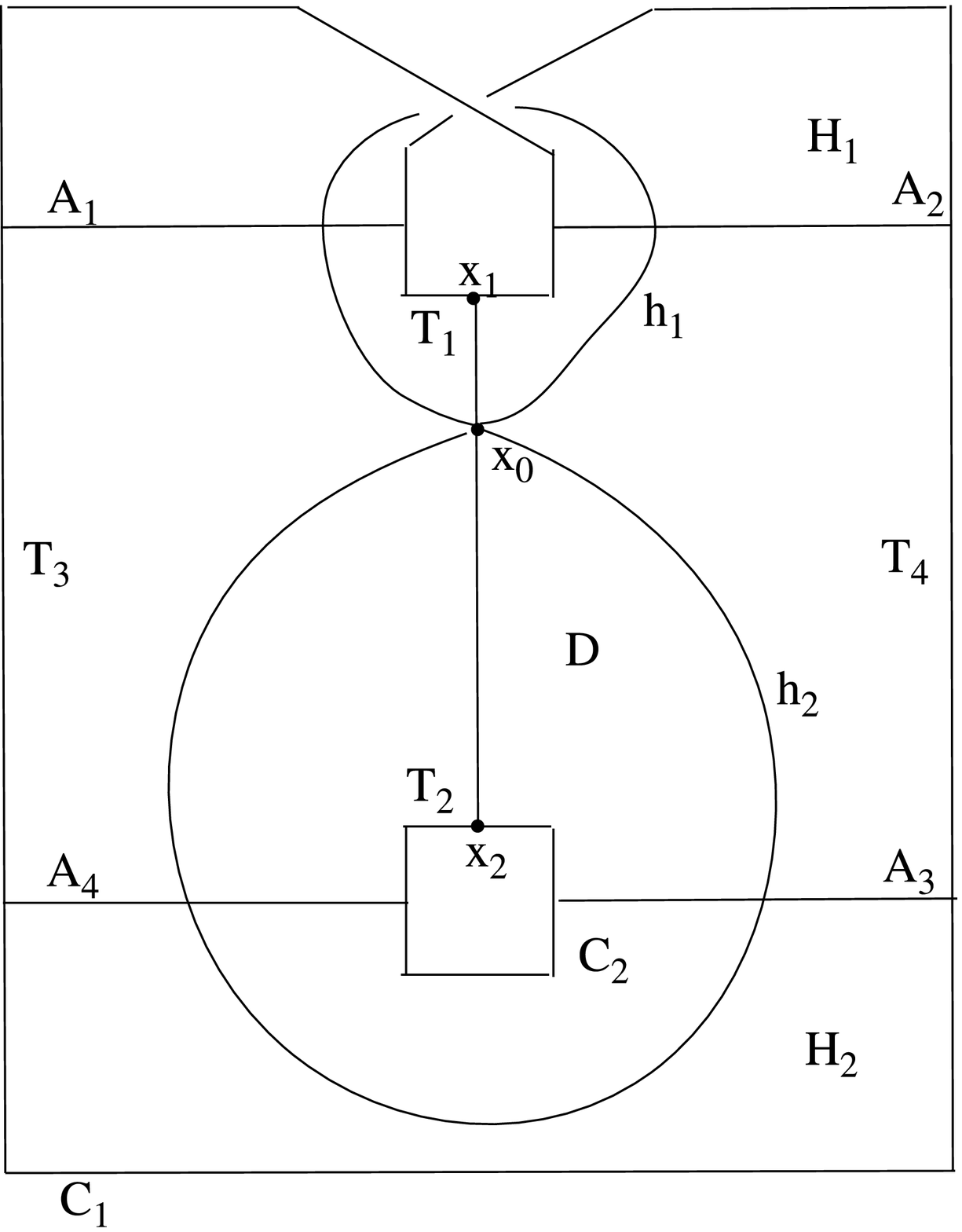}
\caption{} 
\end{figure}

For a map $f \co (X, \partial X) \to (X, \partial X)$, the
restriction of $f$ to the boundary component $C_i$ will always be
called $f_i$ and we will use $C_{i^\#}$ to denote the component of
$\partial X$ containing $f(C_i)$.  Thus, corresponding to the map
$f$,  for $i = 1, 2$ we have the maps of simple closed curves $f_i
\co C_i \to C_{i^\#}$.  The absolute value of the degree of a map
between simple closed curves is independent of the orientations of
the curves, so it can be used without specifying orientations.  For
that reason, by the {\it degree} of a map between simple closed
curves we will mean the absolute value.  In particular, we denote the
degrees of the maps $f_i$ by $|d_i|$.

In order to analyze the boundary-preserving maps of $X$, we will need
the classification of the simple closed curves in the interior of $X$
that are transverse to $A$ and have minimal intersection with it, up
to isotopy.  That is the purpose of the first two lemmas.

\begin{lemma}
Let $K$ be a simple closed curve in the interior of
$X$ that is transverse to $A$ and such that the number of points in
$K \cap A$ is minimal within the isotopy class of $K$.  Then each of
$K \cap H_1$ and $K \cap H_2$ has at most two components.
\end{lemma}

\proof Write $K \cap H_1$ as a union of its components as $K
\cap H_1 = J_1 \cup \dots \cup J_{r_1}$.  Each component $J_j$ is an
arc with its endpoints in $A$.  If there is any components of $K
\cap H_1$ with both endpoints in the same $A _i$, take the innermost
such component, then it bounds a disc in $H_1 - K$ that can be used
to reduce $K\cap A$ by two points, contrary to the assumption of
minimal intersection.  Thus we may assume that each component $J_j$
has endpoints $p_{1j} \in A_1$ and $p_{2j} \in A_2$.  Orienting the 
boundary
of $D$ in a clockwise direction induces an ordering on the $A_i$.
Ordering the $J_j$ so that $p_{1j} < p_{1,j+1}$ in the ordering of
$A_1$ then, since the handle $H_1$ is twisted, $p_{2j} < p_{2,j+1}$
for all $j$.  Let $L_1 \subset K \cap (D \cup H_2)$ be
the component containing $p_{11}$ and let $x$ be the other endpoint of 
$L_1$.
If $x = p_{21}$, then $J_1 \cup L_1$ is a simple closed curve so $
r_1 = 1$.   Now suppose $x \ne p_{21}$.  If $x \ne p_{2{r_1}}$ then
$p_{2{r_1}}$ and the subarc of $A_1$
between $p_{11}$ and $p_{1{r_1}}$ are in different components of $(D 
\cup
H_2) - L_1$ (which is disconnected because both endpoints of the arc
$L_1$ are in the same component of the boundary of the annulus $D \cup 
H_2$).
Since there could not then be an component of $K \cap (D \cup H_2)$
with endpoint $p_{2r_1}$, we conclude that the endpoints of $L_1$ are
$p_{11}$ and $p_{2r_1}$.  A symmetric argument shows that if $L_2$ is
the component of $K \cap (D \cup H_2)$ containing $p_{21}$ then the
other endpoint is $p_{1{r_1}}$.  Since $J_1 \cup L_1 \cup J_2 \cup L_2$ 
is
a simple closed curve, we conclude that $r_1 \le 2$.

Writing $K \cap H_2 = J'_1 \cup \dots \cup J'_{r_2}$, the arcs $J'_j$
have endpoints $p_{3j} \in A_3$ and $p_{4j} \in A_4$.  Ordering the
$J'_j$ so that $p_{3j} < p_{3,j+1}$ in the ordering of $A_3$ induced
by the orientation of the boundary of $D$ then, since the handle $H_2$
is not twisted, $p_{4j} > p_{4,j+1}$ for all $j$.  Letting $L'_1$
denote the component of $K \cap (D \cup H_1)$ containing $p_{31}$, if
the other endpoint is $p_{41}$, then $J'_1 \cup L'_1$ is a simple
closed curve and $r_2 = 1$.   Thus we suppose that endpoint is not
$p_{41}$.  If $L'_1 \subset D$, then $p_{41}$ and the subarc of $A_3$
between $p_{31}$ and $p_{3{r_2}}$ are in different components of $D -
L'_1$. Since there could not then be a component of $K \cap (D \cup
H_1)$  with endpoint $p_{41}$, we conclude that $L'_1 \cap H_1$ is an
arc.  A symmetric argument shows that if $L'_{r_2} \subset D$ then no
component of $K \cap (D \cup H_1)$ could have endpoint $p_{4r_2}$ and
so $L'_{r_2} \cap H_1$ is also an arc.  Since we showed in the first
part of the proof that there are at most two components of $K \cap
H_1$, we conclude that $r_2 \le 2$.
\endproof

Let $K$ be a simple closed curve in the interior of
$X$ that is transverse to $A$.  If $K \cap H_j$ has $r_j$ components
for $j = 1, 2$, then we will say that $K$ is {\it of type} $(r_1, r_2)$.
A simple closed curve in $X$ is called {\it inessential} if it bounds
a disc and {\it essential} otherwise.

\begin{lemma} Let $K$ be an essential simple closed curve in the 
interior of
$X$ that is transverse to $A$ and such that the number of points in
$K \cap A$ is minimal within the isotopy class of $K$.  Then the type
$(r_1, r_2)$ of $K$ is one of the following: $(1, 0), (2, 0), (0, 1),
(1, 1), (2, 1)$ or $(2, 2)$,
\end{lemma}

\proof By Lemma 1, we know that the possible
values of the
$r_j$ are $0, 1$ and $2$.  If $(r_1, r_2) = (0, 0)$ then $K$ would 
bound a
disc, in $D$, which is excluded since $K$ is essential.  We
claim that $K$ cannot be of type $(0, 2)$ or $(1, 2)$.  We assume $r
= 2$ and we establish a contradiction as follows.  Let $J'_1$ be the
component of $K$ in $H_2$ with endpoints $p_{j1}$ for $j = 3, 4$ and
let $L'_1$ be the component in $D \cup H_1$ with endpoint $p_{41}$,
then the other endpoint must be $p_{32}$ (because if it were $p_{31}$
then $J'_1 \cup L'_1$ would be a simple closed curve and if it were
$p_{42}$ that would contradict minimality because either $L'_1$ or
$L'_2$ could be eliminated).  If $L'_1 \subset D$, then $p_{31}$ and
$p_{42}$ are in different components of $(D \cup H_1) - L'_1$ and so
the case $(r_1, r_2) = (0, 2)$ cannot occur.  If $L'_1$ intersects
$H_1$, then, since $p_{11}$ and $p_{21}$ are in different components
of $D - L'_1$, the case $(r_1, r_2) = (1, 2)$ also cannot occur.
The six cases remaining are shown in Figure 2.
\endproof

\begin{figure}[ht!]
\centering
\includegraphics[scale=0.4]{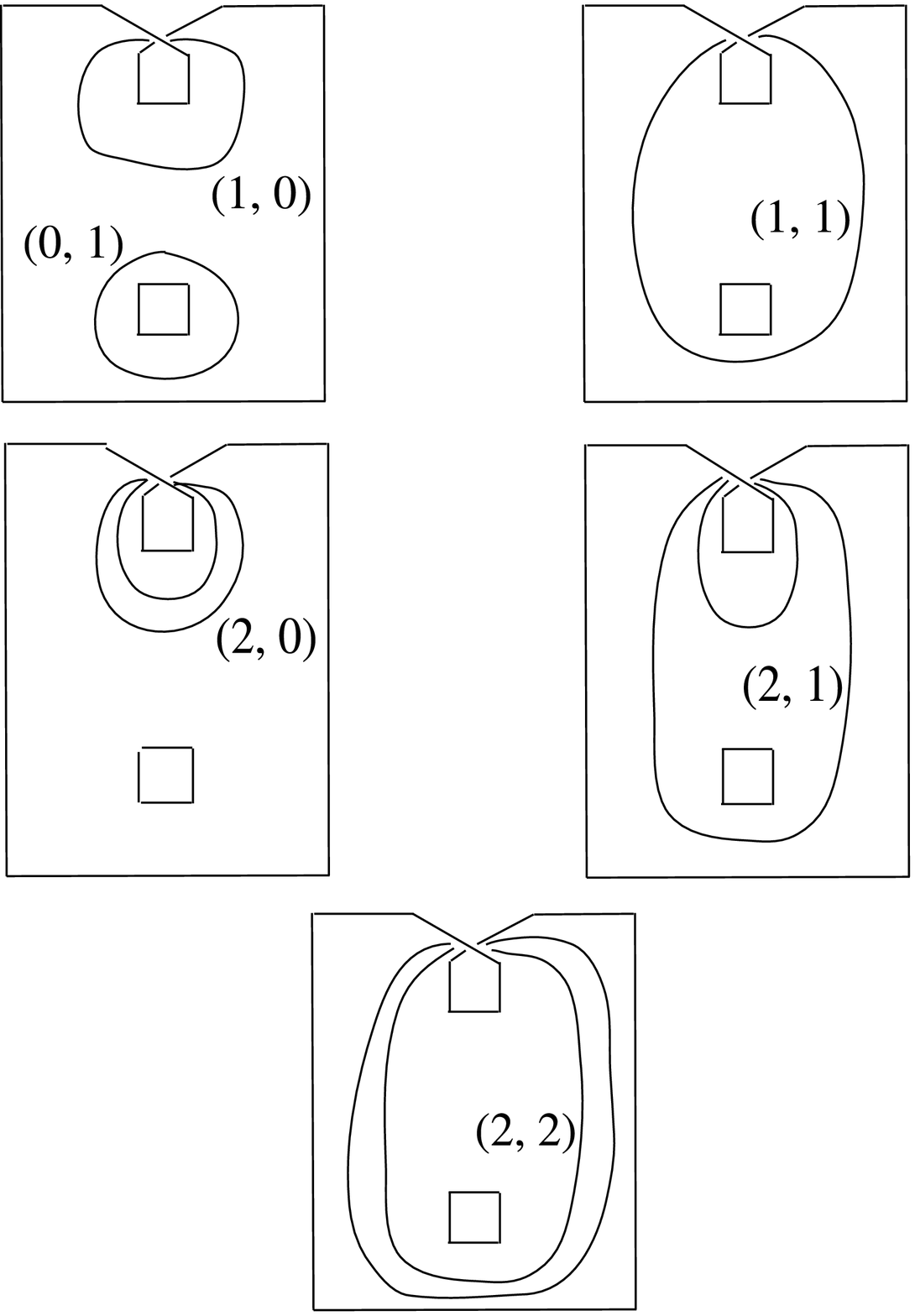}
\caption{}
\end{figure}

We will now describe a convenient form for a map 
$f \co (X, \partial X) \to  (X, \partial X)$ and then prove that we can always obtain it through a boundary-preserving homotopy.  
An arc in $X$ is {\it properly embedded} if it intersects $\partial X$ in its endpoints and at no
other point.  A compact $1$-dimensional submanifold of $X$ is said to be {\it properly embedded} if it is a union of simple closed curves in the interior of $X$ and properly embedded arcs.  
A properly embedded arc $K$ in $X$ is {\it
inessential} if there is an arc $L$ in $\partial X$ connecting the
endpoints of $K$ such that the simple closed curve $K \cup L$ 
bounds a disc in $X$.  
Otherwise, the properly embedded arc $K$ is {\it essential}.

\medskip
{\bf Definition}\qua A map $f \co (X, \partial X) \to  (X, \partial
X)$ is said to be {\it in reduced form} if it has the following
properties: 

(1)\qua for $i = 1, 2$, the map $f_i \co C_i \to C_{i^\#}$
is the constant map to $x_{i^\#}$ if $|d_i| = 0$ and a $|d_i|$-to-one
map otherwise,  

(2)\qua $f^{-1}(A)$ is a properly embedded
$1$-dimensional submanifold transverse to $A$ such that each
component of $f^{-1}(A)$ has minimal intersection with $A$ within its
proper isotopy class, 

(3)\qua all the components of $f^{-1}(A)$ are
essential and 

(4)\qua if $\lambda$ is an arc in $A$ with endpoints in
$f^{-1}(A)$ and interior disjoint from $f^{-1}(A)$, then the
endpoints of $\lambda$ are mapped to distinct components of $A$.

\begin{lemma} Given a map $f \co (X, \partial X) \to  (X, \partial
X)$, there is a map $g$ in reduced form that is homotopic, as a map
of pairs, to $f$.
\end{lemma}

\proof The restriction of $f$ to $\partial X$ can be made to
satisfy condition (1) by the homotopy extension theorem.
Without changing the map on the boundary, and thus condition (1) is
retained, the map can then be modified, using the techniques of
\cite{k3}, so that the inverse image of $A$ is a properly embedded
$1$-submanifold whose components make minimal intersection, within
their proper homotopy class, with $A$.  In this way, (2) is also 
satisfied.
For condition (3), components that are inessential simple closed curves
can be eliminated by an ``innermost simple closed curve" argument and,
similarly, inessential properly embedded arcs can be eliminated by an
``outermost arc" argument.  To simplify the notation, let us suppose
that $f$ already satisfies conditions (1), (2) and (3) of the
definition of reduced form, then it remains to homotope $f$ so that,
if $\lambda$ is an arc in $A$ with endpoints in
$f^{-1}(A)$ and interior disjoint from $f^{-1}(A)$, then the
endpoints of $\lambda$ are mapped to distinct components of $A$.   To
accomplish this, suppose the endpoints of such an arc $\lambda$ are
mapped to the same component of $A$.  Since the interior of $\lambda$
does not contain points of  $f^{-1}(A)$, the path $f(\lambda)$ either
lies entirely in $D$ or entirely in one of the handles $H_j$.  If
$f(\lambda)$ lies in $D$ then, since $f^{-1}(A)$ is transverse to
$A$, there are neighborhoods of the endpoints of $\lambda$ in $A$
whose image lie in some handle $H_j$.  We can homotope $f$ so that
the image of a neighborhood of $\lambda$ in $A$ lies entirely in that
handle.  Similarly, if $f(\lambda)$ lies in some handle $H_j$, then
we can homotope $f$ so that the neighborhood of $\lambda$ is mapped
to $D$.  Repeating this procedure a finite number of times we obtain
a map in reduced form homotopic to $f$.
\endproof

\begin{theorem}
Let $f \co (X, \partial X) \to  (X, \partial X)$ be a map in
reduced form.  Then $f$ may be homotoped to a map in reduced form
with the property that every simple closed curve in the inverse image
of $A$ is of type $(0, 1)$.
\end{theorem}

\proof  We
claim that all the types of simple closed curves listed in Lemma 2
except $(0, 1)$ may be excluded as a component of the inverse image
of $A$ for some map
in reduced form homotopic to the given map.
Suppose $K$ is a component of $f^{-1}(A)$ of type $(1, 0)$ or type
$(1, 1)$.  Then there is a closed neighborhood $N$ of $K$ such that
$N \cap f^{-1}(A) = K$ and
$N - K$ is connected.  Since $f(N - K)$ is a connected subset of
$X - A$, it lies in one of the three components of $X - A$, each of
which is contractible.  We denote the component of $X - A$ that
contains $f(N - K)$ by $B$.  Since $f$ maps the boundary $\partial
N$ of $N$ into the contractible set $B$, we can extend $f|\partial N$
to a map $g \co N \to B$.  Extend $g$ to all of $X$ by letting it
equal $f$ outside of $N$.  Then $g$ is a reduced form of $f$ with
$g^{-1}(A) = f^{-1}(A) - K$.  Repeating this procedure a finite
number of times, we may assume we have a reduced form for $f$ (which
we still call $f$) in which no component of $f^{-1}(A)$ is a simple
closed curve of type $(1, 0)$ or type $(1, 1)$.  Since $f$ is in
reduced form, if $\lambda$ is an arc in $A$ with endpoints in
$f^{-1}(A)$ and
interior disjoint from $f^{-1}(A)$, then the endpoints of $\lambda$
must be mapped to distinct components of $A$.  We claim that this
property implies that
simple closed curves of types $(2, 0),
(2, 1)$ and $(2, 2)$ cannot be components of $f^{-1}(A)$.
To prove it, let $K \subseteq f^{-1}(A)$ be an innermost simple
closed curve of type
$(2, n)$ for any $n = 0, 1, 2$, then $f(K)$ is connected
so it lies entirely inside one of the $A_j$.  Taking the two points
of intersection of $K$ with $A_1$, the arc
$\lambda \subseteq A_1$ between them would have no points of the
inverse image of $A$ in its interior because there are now no simple
closed curves of types $(1, 0)$ or $(1, 1)$.  But the endpoints of
$\lambda$ map to the same
component of $A$, contrary
to the definition of reduced form.
Since we have eliminated from the inverse image of $A$ all
simple closed curves of types $(0, 1), (1, 1), (2, 0), (2, 1)$ and
$(2, 2)$, we conclude that, if there are any simple
closed curves, they must be of type $(0, 1)$.
\endproof

In order to complete the analysis of the components of $f^{-1}(A)$
for $f$ in reduced form, we turn now to the isotopy classification of
properly embedded arcs in $X$.  Recall that the
components of $\partial D - (H_1 \cup H_2)$ are $T_1 \cup T_2 \cup T_3
\cup T_4$ where $T_2 \subset C_2$, as in Figure 1.  If a properly
embedded arc $K$ is isotopic to a line segment in $D$ with endpoints
in $T_i$ and $T_j$, then we will call $K$ an arc {\it of type} $[i, j]$.

We begin with the
classification of arcs $K$ whose endpoints lie in different
components of  $\partial X$.

\begin{theorem}
Let $K$ be an arc in $X$ that has minimal intersection with $A$ in
its isotopy class and whose endpoints lie in different
components of $\partial X$.
Then $K$ is of type $[1, 2]$ or of type $[2, 3]$ and $K$ lies entirely 
in $D$.
\end{theorem}

\proof We can isotope $K$ so that one endpoint lies in $T_1
\cup T_3 \cup T_4$
and the other in $T_2$ and $K$ has minimal
intersection with $A$ within its isotopy class.  It will be
convenient to orient $K$ so that we view $K$ as
ending in $T_2$.  We will first show that the first intersection of
$K$ with $H_1 \cup H_2$ cannot be with $H_1$.  The intersections of
$K$ with $\partial D$ cannot begin $T_1A_1A_2$ because $K$ could be
isotoped to begin at $T_4$, so the
intersection with $A$ is not minimal.  In the same way, the
intersection  $T_1A_2A_1$ can be isotoped to begin at $T_3$.  If the
intersections began with $T_3A_1A_2$ we can isotope to begin at $T_1$
whereas, if it began with $T_3A_2A_1$, the arc
could not reach $T_2$.  Checking the two
corresponding cases for $K$ beginning in $T_4$, we conclude that if
$K$ intersects  $H_1 \cup H_2$, it must first intersect $H_2$.  If
$K$ starts in $T_1$, intersects $H_2$ and then $H_1$, we can assume
without loss of generality that the intersections are
$T_1A_3A_4A_1A_2$, but then $K$ cannot end in $T_2$.  On the other
hand if, after
intersecting $H_2$, the arc $K$ does not next intersect $H_1$, then
it must have its last intersection with $H_2$.  But it is clear that,
if the last intersection is with $H_2$, we can isotope $K$ to reduce
the intersection with $A$, so that intersection was not minimal.
Thus, if $K$ starts at $T_1$ and has minimal intersection with $A$,
then it does not intersect $A$ at all.  If $K$ starts at $T_3$ and
intersects $A$ in the order $T_3A_3A_4$, it also
would have its last intersection with $H_2$ and thus not have minimal
intersection with $A$.  The intersection with $A$ is also not minimal
if $K$ begins $T_3A_4A_3$ because the first intersection can be
reduced by isotoping $K$ to begin in $T_4$.
Symmetric arguments when $K$ begins in $T_4$ lead us to conclude
that, no matter where $K$ begins, it cannot intersect $H_1 \cup H_2$
and so $K$ lies entirely in $D$.  Therefore $K$ is isotopic to a line
segment ending in $T_2$ and beginning in either $T_1, T_3$ or $T_4$.
Since the line segment between $T_2$ and $T_4$ is
isotopic to the line segment between $T_2$ and $T_3$, there are just
the two possibilities given.
\endproof

\begin{figure}[ht!]
\centering
\includegraphics[scale=0.4]{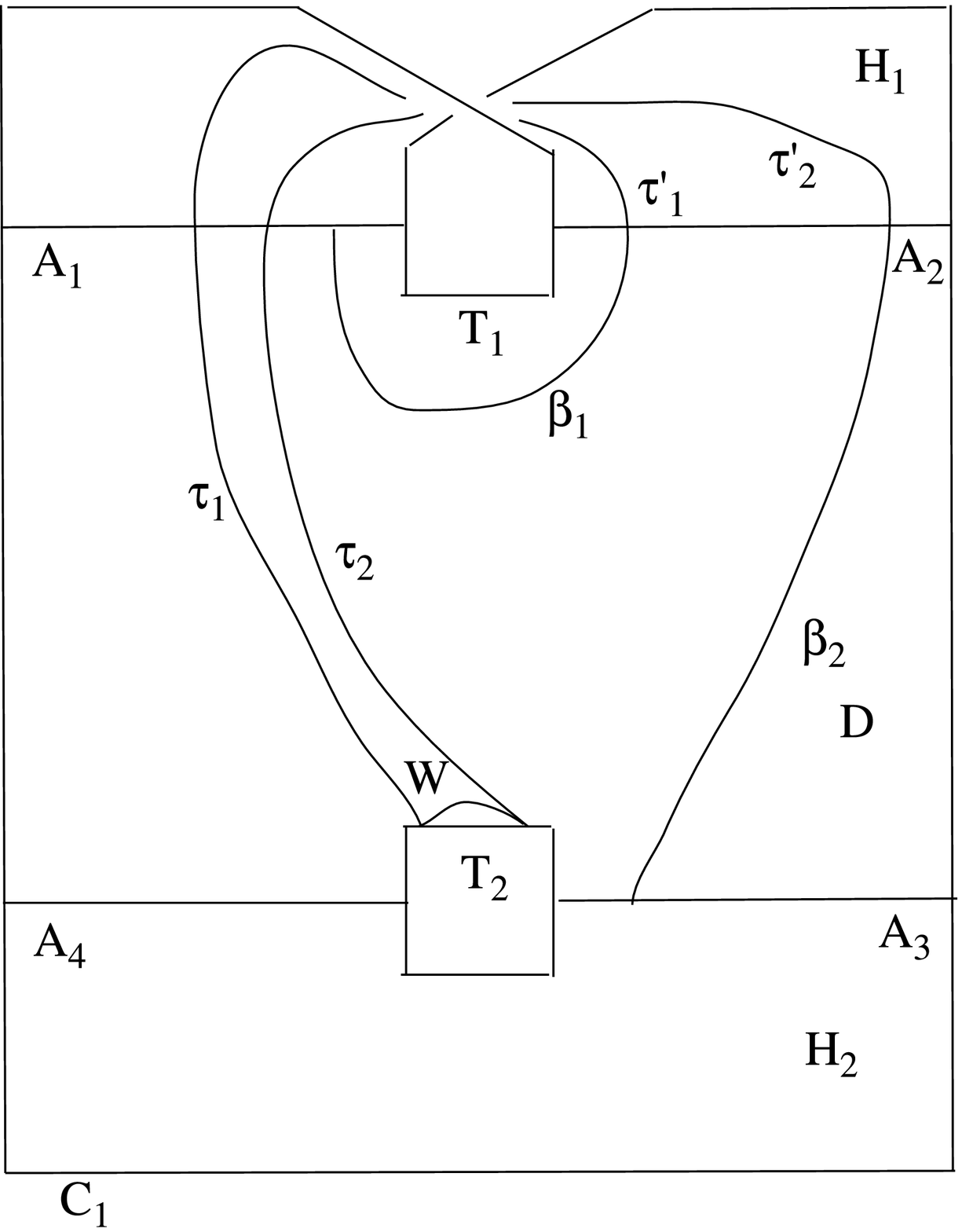}
\caption{} 
\end{figure}

\begin{lemma}
Let $K$ be an arc in $X$, both of whose endpoints lie in $C_2$,
such that the number of points in $K \cap A$ is minimal
in the isotopy class of $K$. Let
$W$ be an arc in $D$ connecting the endpoints of $K$ and otherwise
disjoint from $K$.  Then the number of points in $(K \cup W) \cap A$
is minimal in the isotopy class of the simple closed curve $K \cup W$.
\end{lemma}

\proof
Let $\tau_1$ and $\tau_2$ be the two components of $K \cap D$ that
intersect $C_2$ (see Figure 3).  Assuming, without loss of
generality, that $\tau_1$ intersects $A_1$.  We note that the
minimality condition on $K \cap A$ also assures us that $\tau_2$ cannot
intersect $A_3$ or $A_4$.   Suppose $\tau_2$ also intersects $A_1$.
Consider those portions of $K$, which we denote by $\tau'_i$ for $i =
1, 2$, that
start at the endpoints in $A_1$ of the respective $\tau_i$ and that
intersect the $A_i$ in the same order, and that are maximal with 
respect to
this property.  Figure 3 illustrates the case in which in the
$\tau'_i$ intersect in the order $A_1A_2$.  The other
possible case, in which the order is $A_1A_2A_3A_4A_1A_2$, can be
described by a similar figure.  Let
$\beta_i$ denote the components of $K - A$ whose closures intersect
the $\tau'_i$.  Since, by definition of the $\tau'_i$, the arcs
$\beta_1$ and $\beta_2$ cannot have their endpoints in the same
$A_i$, one must have an endpoint in $A_1$ and the other an endpoint
in $A_3$, as we illustrated in Figure 3.  Relative to those
endpoints, we can isotope the portion of $K \cup W$ consisting of
$W \cup \tau_1 \cup \tau_2 \cup \tau'_1 \cup \tau'_2 \cup \beta_1 \cup
\beta_2$ to an arc in $D$ with endpoints in $A_1$ and $A_3$.  Thus,
if we isotope $K \cup W$ so that its intersection with $A$ is
minimal, we obtain a simple closed curve that contains such an arc.
However, as we can see in Figure 2, none of
the simple closed curves that have minimal intersection with $A$ can
contain an arc in $D$ connecting $A_1$ and $A_3$.  Thus, assuming
that $\tau_2$ intersects $A_1$ has lead to a contradiction and we
conclude that $\tau_2$ intersects $A_2$, which implies that the
number of points in $(K \cup W) \cap A$ is minimal in the isotopy
class of the simple closed curve $K \cup W$.
\endproof

The punctured M\"obius band $X$ can be constructed
from a surface $\widetilde X$ that consists of the two-dimensional
sphere with four symmetrically-placed open discs deleted, by
identifying antipodal points.  Consider a 90-degree
rotation of $\widetilde X$ that is equivariant with respect to the
antipodal action.  Then, for $\phi \co X \to X$ the corresponding
homeomorphism, we have $\phi(C_2) = C_1$ and $\phi(C_1) = C_2$.
Moreover, we may assume that $\phi(x_1) = x_2$ then, since $\phi^2$
is the identity map, we have $\phi(x_2) = x_1$ as well.

\begin{lemma}
Let $K$ be an arc in $X$ both of whose endpoints lie in the same
component of $\partial X$.  If both endpoints are in $C_2$, then $K$
is isotopic to an arc obtained by removing an arc from a
simple closed curve of type $(1, 0), (2, 0)$ or $(2, 1)$ and connecting 
the
endpoints of the arc that remains to $C_2$ by parallel arcs.
Such arcs are said to be of types $\{1, 0\}, \{2, 0\}$ and $\{2, 1\}$,
respectively.  If both
endpoints are in $C_1$, then either $K$ is of type $[1, 3]$ or
of type $[3, 4]$, or it is isotopic to an arc with both endpoints in
$T_1$ that passes once through $H_2$ and does not intersect $H_1$.
Such an arc is said to be of type $\{0, 1\}$.
\end{lemma}

\proof Suppose both endpoints of $K$ are in $C_2$.  Isotope $K$
so that the number of points in $K \cap A$ is minimal in the isotopy
class of $K$.  Let $W$ be an arc in $D$ connecting the endpoints of
$K$ and otherwise disjoint from $K$.  By Lemma 4, the number of
points in $(K \cup W) \cap A$ is minimal in the isotopy class of the
simple closed curve $K \cup W$ and therefore $K \cup W$ is of one of
the six types of simple closed curves described by Lemma 2.  Thus,
once the arc $K$ has been isotoped so that the number of points in $K
\cap A$ is minimal in the isotopy class, it must have the same
intersections with $H_1$ and $H_2$ as in that lemma.  We can
describe $K$ itself by deleting an arc
in $D$ from a simple closed curve of
Lemma 2 and connecting the endpoints of what remains with $C_2$
by parallel arcs.  It is easy to check that, if we carry out
this construction with the types $(0, 1), (1, 1)$ or
$(2, 2)$, then $K \cap A$ is not minimal in the isotopy class of $K$.
That leaves the three types  $(1, 0), (2, 0)$ and $(2,
1)$; see Figure 4.

Now suppose both endpoints of $K$ are in $C_1$.  Then the endpoints
of $\phi(K)$ are in $C_2$ and so $\phi(K)$ is isotopic to some arc $K'$
that is one of the three types just described.  Applying
$\phi = \phi^{-1}$ once more, we have $K$ isotopic to $\phi(K')$ and
thus there are three isotopy types of arcs with both endpoints in
$C_1$.  To show that those three isotopy types are the ones described
in the statement of the theorem, it suffices to prove that no two
such arcs are boundary-preserving isotopic.   Proper arcs are relative
$1$-cycles that represent classes of $H_1(X, \partial X)$.
If two properly embedded arcs are boundary-preserving isotopic, then
the corresponding cycles are homologous.  However, the
arcs of type $[1, 3], [3, 4]$ and $\{0, 1\}$ represent three
different elements of
$H_1(X, \partial X)$.
\endproof

In the proof of Lemma 5, we used the homeomorphism $\phi$, that
interchanges the components of $\partial X$, in order to apply our
analysis of the isotopy  classes of arcs with endpoints in $C_2$ to
arcs whose endpoints lie in $C_1$.  We will also use $\phi$ so that
information about minimizing the number of fixed points through a
homotopy can be extended to additional cases by means of the
following, easily verified, result.

\begin{lemma} Let $f \co (X, \partial X) \to  (X, \partial X)$ be
a map and suppose that $g \co (X, \partial X)$ $\to  (X, \partial
X)$ is homotopic, as a map of pairs, to $\phi f \phi$, then $\phi g
\phi$ is homotopic to $f$ and the set of fixed points of $\phi g
\phi$ is homeomophic to the set of fixed points of $g$.
\qed \end{lemma}

If a map $f \co (X, \partial X) \to  (X, \partial X)$ is in
reduced form, then for $i = 1, 2$, the 
map $f_i \co C_i \to C_{i^\#}$ is either the
constant map at $x_{i^\#} \notin A$ if $|d_i| = 0$ or, otherwise,
$f_i$ is a $|d_i|$-to-one map. If $|d_i| = 0$, then $f_i^{-1}(A) \cap 
C_i$
is the empty set.

\medskip
{\bf Definition}\qua Suppose that $|d_i| \ne 0$ then, for each point
of the finite set $f_i^{-1}(A) \cap C_i$ we record the integer $1, 2,
3$ or $4$ corresponding to the component of $A$ to which it is sent by
$f$.  As we go once around the simple closed curve $C_i$, we obtain an
ordered set made up of the integers 1 through 4 which we call the {\it
pattern} of $f_i$.  Of course the pattern depends on where on $C_i$ we
begin our circuit of the simple closed curve and in which direction we
go.

\begin{lemma}
Given a map $f \co (X, \partial X) \to  (X, \partial X)$ in
reduced form, if the degree $|d_i|$ of $f_i \co C_i \to C_{i^\#}$  is 
non-zero, then we may choose the starting point and direction so that 
the pattern of $f_i$ is
$(121234)^{|d_i|}$ if $i^\# = 1$ and $(34)^{|d_i|}$ if $i^\# = 2$.
\end{lemma}

\proof Since $C_1$ intersects $A$ in six points then,
if $i^\# = 1$, it must be that $f_i^{-1}(A) \cap C_i$ consists of 
$6|d_i|$
points because $f$ is in reduced form.  Moreover, by starting at a
suitable place in $C_i$ and choosing the correct direction, the
points of $f_i^{-1}(A) \cap C_i$
that we encounter will be first mapped to $A_1$, then to $A_2$, next
to $A_1$ again and
so on to produce the pattern $(121234)^{|d_i|}$ (see Figure 1).  On
the other hand, $C_2$ intersects $A$ only at one point each of $A_3$
and $A_4$.  Therefore, if  $i^\# = 2$ then, choosing the correct
starting point and direction, the pattern in this case will be
$(34)^{|d_i|}$.
\endproof

\begin{theorem}
Let $f \co (X, \partial X) \to  (X, \partial X)$ be a map in
reduced form such that $f_1$ is essential, then every proper arc in
$f^{-1}(A)$ is one of the following types:
$[1, 2], [1, 3], [2, 3]$ or $\{1, 0\}$.
\end{theorem}

\proof  We first show that there are no arcs of type $[3, 4]$
in $f^{-1}(A)$.  Suppose that there are and let $K$ by the
``innermost" such arc.  That is, letting $L$ be the component of $X -
K$ that contains $C_2$, then no other arc of that type is in $L$.
Arcs that are components of $f^{-1}(A)$ that lie in $L$ must be of
type $[2, 3]$.  Let $f(K) = A_j$.  If there were no arcs of type $[2,
3]$ in $L$, the pattern of $f_1$ would contain $jj$, contrary to the
Lemma 7.  Thus there must be type $[2, 3]$ arcs in $L$ and these are
the only components of $f^{-1}(A)$ that intersect $C_2$.   Therefore
the number of type $[2, 3]$ arcs is either $6|d_2|$ or $2|d_2|$,
depending on whether $C_2$ is mapped to $C_1$ or to $C_2$.  But that
is an {\it even} number in either case whereas, by Lemma 7, the
pattern of $f_1$ is either $(121234)^{d_1}$ or $(34)^{|d_1|}$.  A
contradiction occurs because, in either possible pattern for $f_1$,
the number of integers between any repetition of an integer $j$ is
{\it odd}.  Thus there are no arcs of type $[3, 4]$.  Basically the
same argument demonstrates that arcs of $\{0, 1\}$ cannot be
components of $f^{-1}(A)$, with arcs of type $[1, 2]$ taking the role
of the arcs of type $[2, 3]$ in the previous argument.

We next eliminate arcs of $\{2, 0\}$ as possible components of
$f^{-1}(A)$.  Suppose there are arcs of this type and let $K$ denote
the ``innermost" such component.  That is, there is an arc $\alpha$
in $C_2$ connecting the endpoints of $K$, and with its interior in
the component of $X - K$ that does not contain $C_1$, such that no
other arc of type $\{0, 2\}$ intersects $\alpha$.  Let $f(K) = A_j$.
If there were no other component of $f^{-1}(A)$ with endpoints in the
interior of $\alpha$, then the pattern of $f_2$ would contain $jj$,
contrary to Lemma 7.  However, a component of $f^{-1}(A)$ with
endpoints in $\alpha$ would be of type $\{1, 0\}$ and the innermost
of these would, in the same way, produce a pattern for $f_2$ contrary
to Lemma 7.  Thus there are no arcs of type $\{2, 0\}$ that are
components of $f^{-1}(A)$.

The proof that there are no components of $f^{-1}(A)$ of type $\{2,
1\}$ is similar to the argument for type $\{2, 0\}$.  This time let
$K$ be the ``outermost" component of type $\{2, 1\}$.  This means
that there is an arc $\beta$ in $C_2$ connecting the endpoints of
$K$, and with its interior in the component of $X- K$ that does not
contain $C_1$, such that no other arc of type $\{2, 1\}$ has its
endpoints in $\beta$.  If there were no components of $f^{-1}(A)$
intersecting the interior of $\beta$, the pattern of $f_2$ would
again contain $jj$ where $f(K) = A_j$, contrary to Lemma 7.  If
there were such components, however, they would be of type $\{1, 0\}$
and the ``outermost" of these would lead to the same contradiction.
Thus there are no components of $f^{-1}(A)$ that are arcs of type
$\{2, 1\}$.  By Theorem 5 and Lemma 5, the four types listed in
the statement of the theorem are the only possibilities remaining
(see Figure 4).
\endproof

\begin{figure}[ht!]
\centering
\includegraphics[scale=0.4]{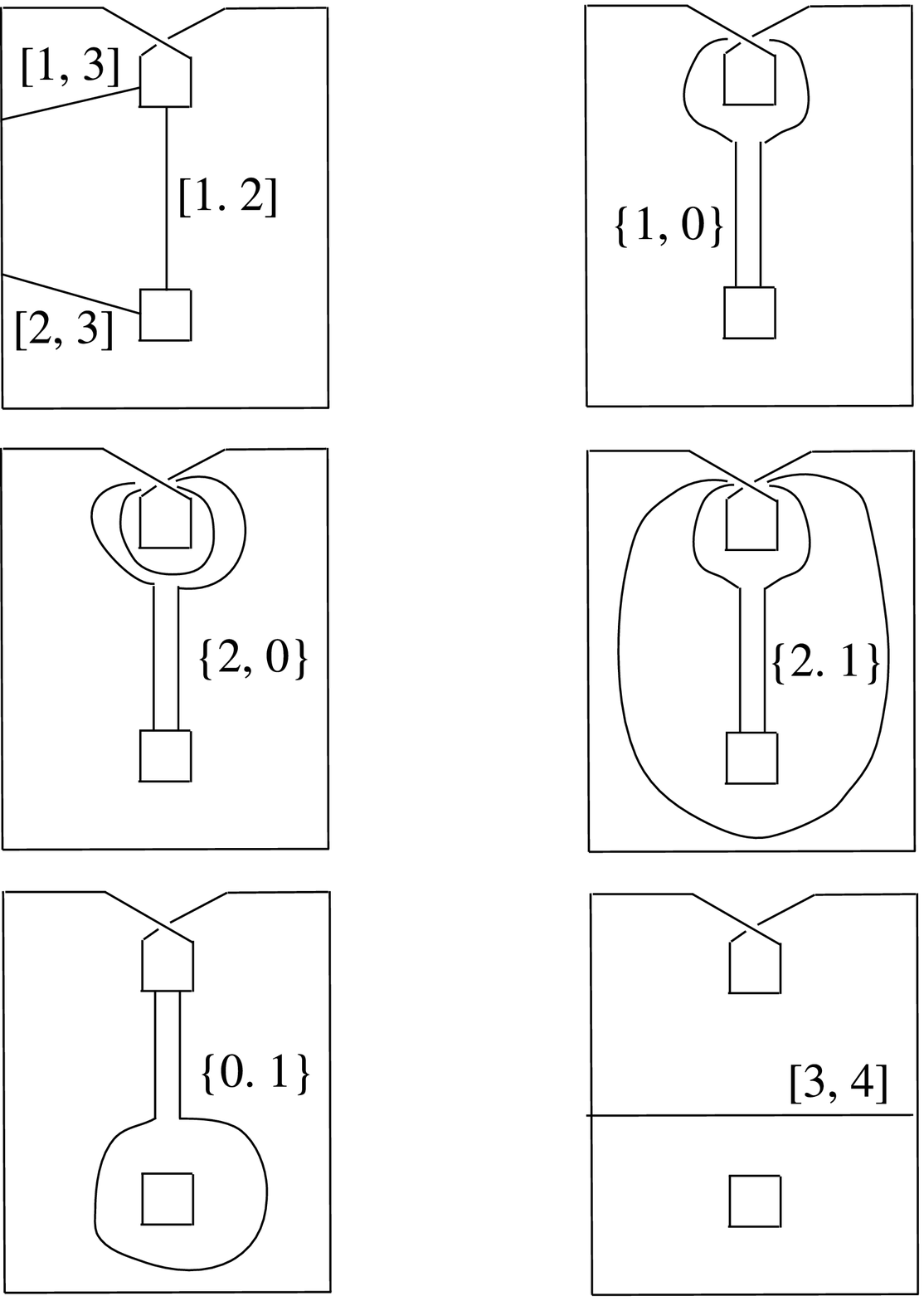}
\caption{}
\end{figure}

\section{Boundary inessential maps}

{\bf Definition} Represent the points of the unit circle $S^1$ by
$e^{i\theta}$ for $0 \le \theta \le 2\pi$.
Let $R$ be an annulus embedded in $X$ and $\alpha$
an arc in $R$ such that there is a homeomorphism $\psi \co [0, 1]
\times S^1 \to R$ with 
$\psi([0, 1] \times \{e^{i0} = 1\}) = \alpha$.  Given a
map $g \co \alpha \to X$, the {\it annulus extension} of $g$ is
the map $g \co R \to X$ defined by $g(\psi(t, e^{i\theta})) =
g(\psi(t, 1))$ for all $0 \le t \le 1$ and $0 \le \theta \le 2\pi$.

\begin{theorem}
If $f \co (X, \partial X) \to  (X, \partial X)$ is boundary
inessential, then $MF_\partial[f] = N_\partial(f)$.
\end{theorem}

\proof By Lemma 3, we may assume that $f$ is in reduced form.  
Therefore
the restriction of $f$ to each
component $C_j$ is a constant map to one of the points $x_1$ or
$x_2$; we write $f(C_j) = x_{j^\#}$.  Consequently, the
components of $f^{-1}(A)$ must be simple closed curves that, by
Theorem 4, we may assume are of type $(0, 1)$.  See Figure 5.   We
denote the components by $K_1, \dots , K_m$ where $K_{i+1}$ is the
component of $X - K_i$ containing $C_2$ and $K_i$ and $K_{i+1}$ bound
an annulus $R_i$ that contains no other component of $f^{-1}(A)$.  We
will describe a map $g \co (X, \partial X) \to  (X, \partial X)$ 
which is
boundary-preserving homotopic to $f$, is also constant on each
component of $\partial X$, and that has a single fixed point in the
interior of $X$ if $Im_\partial(f) = 2$ and $f(C_1) = C_2$.   In that
case, $g$ has no fixed points on $\partial X$ and, since its
Lefschetz number is nonzero, we have $MF_\partial[f] = N_\partial(f)$.
Otherwise, $g$ will have no fixed points in the interior of $X$ so
obviously $MF_\partial[f] = N_\partial(f)$.

Let $[x_1,x_2]$
denote the line segment and let $h_1$ and $h_2$ be the simple
closed curves shown in Figure 5.

We define a map $g'$ that sends all of $X$
to $G = [x_1, x_2] \cup h_1 \cup h_2$.
Let $R$ be the annulus in $X$
bounded by the simple closed curves $h_2$ and $C_2$.  We map the
closure of $X - R$ to the line segment $[x_0, x_{1^\#}]$ by sending
$h_2$ to $x_0$ and extending.
We will obtain $g'$ on $R$ by defining a map on the line segment $[x_0,
x_2]$ and then extending to all of $R$ by using annulus extensions.
The behavior of the intersection of $[x_0, x_2]$ with each
$R_i$ will determine how $g'$ is defined on $R_i$.  We must send $x_2$
to $x_{2^\#} = f(C_2)$, we let $g'(x'_0) = x_0$ and we extend linearly
to the line segment $[x_2,  x'_0]$.  Now each annulus $R_i$ is mapped
by $f$ either to $D$ or to one of the handles $H_j$ and, for adjacent
annuli $R_i$ and $R_{i+1}$, one must go to $D$ and the other to a
handle.  For each $R_i$ mapped to $D$, choose a point $x_i$ in the
interior of the line segment $R_i \cap [x_0, x'_0]$.  For each such
$i$, the map $g'$ sends the line segment $[x_i, x_{i+2}]$ once around
$h_1$ or once around $h_2$, with $x_i$ and $x_{i+2}$
sent to $x_0$, depending on
whether $f$ maps $R_{i+1}$ to $H_1$ or to $H_2$.
Annulus extension on each $R_i$ completes the definition of $g'$.

\begin{figure}[ht!]
\centering
\includegraphics[scale=0.4]{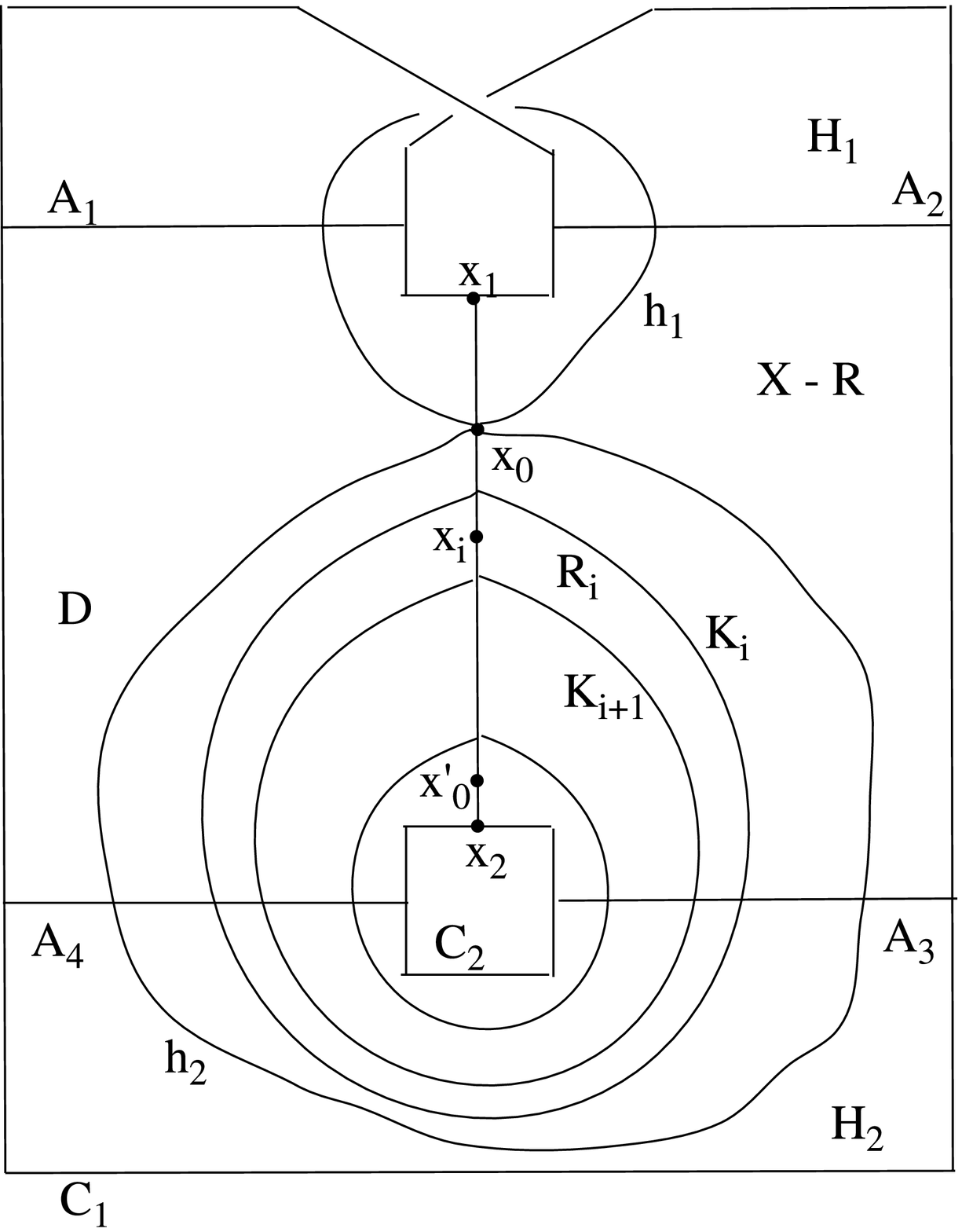}
\caption{} 
\end{figure}

If $f(C_1) = x_1$, then $g'$ is the identity on the line segment $[x_0,
x_1]$ but there are no other fixed points on the interior of $X$.
Let $G' = [x_1, x_2] \cup h'_1 \cup h'_2$ where $h'_i$ is isotopic to
$h_i$ and intersects $[x_1, x_2]$ only at $x_1$.  Let $\rho \co G
\to G'$ map $[x_1, x_2]$ to itself, taking $[x_1, x_0]$ to $x_1$ and
extending linearly over $[x_0, x_2]$, and sending each $h_i$ to $h'_i$
by a homeomorphism.
Note that $\rho g'$ is homotopic to $f$ as a map of the pair
$(X, \partial X)$.   If $f(C_1) = x_1$ and $f(C_2) = x_2$, then $x_1$
and $x_2$ are the only fixed points of $g = \rho g'$.  If
$f(C_1) = f(C_2) = x_1$, then $g = \rho g'$ has a fixed point only
at $x_1$.  Now suppose $f(C_1) = f(C_2) = x_2$ then, by Lemma 6,
the map $g = \phi \rho g'\phi$ is homotopic to $f$ and its only fixed
point is $x_2$.  Finally, if $f(C_1) = x_2$ and $f(C_2) = x_1$, then
we can set $g = g'$ because $g'$ will have the required property, 
namely,
one fixed point in the interior of $X$, at $x_0$. \endproof

\section{Boundary essential maps}

\begin{theorem}
If $f \co (X, \partial X) \to  (X, \partial X)$ is a
boundary essential map such that $Im_\partial(f) = 2$,
then $N_\partial(f) = MF_\partial[f]$.
\end{theorem}

\proof Suppose for now that $f(C_i) = C_i$ for $i = 1, 2$.
Throughout the proof, homotopies will be understood to be relative to 
$\partial X$.  We will show that $f$ maps $X$ onto itself.  Let $Z = 
\partial X \cup A_1 \cup A_3$.  Suppose $f$ is not onto, then we may 
assume there is a point of $X - Z$ that is not in the image of $f$.  By 
retracting, there is a map homotopic to $f$ whose image is $Z$.   
We use \cite{brs}, to 
obtain a map $g \co X \to Z$ such that, for 
$i = 1, 2$, the set $L_i = g^{-1}(x_i)$ is a 
$1$-dimensional subpolyhedron of $X$.
There is a component $\alpha_i$ 
of $L_i$ that is an essential properly embedded arc or an essential 
simple closed curve 
intersecting $C_i$ in a single point or else 
connected to $C_i$ by an arc meeting $C_i$ in a point.   If no such 
$\alpha_i$ existed, there would be a simple closed curve 
in the complement of $L_i$
that is isotopic to $C_i$.  But that 
would imply that $f_i \co C_i \to C_i$ is inessential, contrary to 
the hypotheses.  If $\alpha_1$ and $\alpha_2$ are both essential 
properly embedded arcs, then it follows from Lemma 5 that $\alpha_1 
\cap \alpha_2 \ne \emptyset$ (see Figure 4).  If one or both of the 
$\alpha_i$ is a simple closed curve intersecting $\partial X$ in one 
point  or connected to it by an arc and 
it were true that $\alpha_1 \cap \alpha_2 = \emptyset$, 
then there would be nearby 
properly embedded arcs with the same property, 
contrary to Lemma 5, so we conclude 
that $\alpha_1 \cap \alpha_2 \ne \emptyset$.  But $L_1 \cap L_2 = 
\emptyset$ by definition so the contradiction proves that $f$ maps onto 
$X$.  Every map of $X$ that is homotopic to $f$ through 
boundary-preserving homotopies is boundary-essential and therefore 
onto, as we have shown, so the geometric degree of $f$ is nonzero.  
Consequently, by a theorem of Edmonds (see \cite{sk}, Theorem 1.1 for a 
corrected version of the result), $f$ is homotopic to $qp$ where $p 
\co X \to Y$ is a pinch and $q \co Y \to X$ is an $n$-fold 
branched cover.  Triangulate $X$ so that the branch points are vertices 
and let $v$ be the number of vertices in the triangulation.  
Triangulate $Y$ so that $q$ is a simplicial map, then the number of 
vertices of this triangulation of $Y$ is $nv - b$ for some $b \ge 0$ 
depending on the branching.  Therefore, the Euler characteristic of $Y$ 
is $\chi(Y) = n\chi(X) - b = -n - b < 0$.  If $p$ were a nontrivial 
pinch, then there would be a subsurface $P$ of $X$, not a disc, such 
that $\partial P$ is connected and $Y = X/P$.  But that would imply 
$\chi(Y) = \chi(X) - \chi(P) + 1 =
-\chi(P) \ge 0$ so the pinch must be trivial, that is, $Y = X$ and $p$ 
is the identity map.  But then we have $-1 = -n - b$ which implies that 
$b = 0$ and $n = 1$, that is, $q$ is a homeomorphism.  It follows by
\cite{jg} that $N_\partial(q) = MF_\partial[q]$ and therefore, since 
$f$ is homotopic to $q$, that
$N_\partial(f) = MF_\partial[f]$ when $f(C_i) = C_i$ for $i = 1, 2$.
Now suppose $f(C_1) = C_2$ and $f(C_2) = C_1$.  Let $\phi \co X
\to X$ be the homeomorphism that
interchanges $C_1$ and $C_2$.  The map $\phi f$ has the property that
$\phi f(C_i) = C_i$ for $i = 1, 2$ so, as we have just demonstrated,
$\phi f$ is homotopic to a homeomorphism
$q$.  But then $f$ is homotopic to the homeomorphism
$\phi q$ so we may again apply  \cite{jg} to complete the proof.
\endproof

\begin{theorem}
If $f \co (X, \partial X) \to  (X, \partial X)$ is a
boundary essential map such that $Im_\partial(f) = 1$,
then $N_\partial(f) = MF_\partial[f]$.
\end{theorem}

\proof  We first consider the case that $f(\partial X) = C_2$.
We claim that $f$, assumed to be in reduced form, is homotopic modulo
the boundary to a map $g \co X \to C_2$ which agrees with $f$ on
$\partial X$ and that would establish $N_\partial(f) =
MF_\partial[f]$ in this case.   Since $f$ is an essential map on
$C_2$, there are components of $f^{-1}(A)$ that are arcs with at
least one endpoint in $C_2$.  No simple closed curve of type $(0, 1)$
can lie in the complement of such an arc; see Figure 4.
Thus there are no simple closed curves in
$f^{-1}(A)$ since, by Theorem 4, we need consider only those of type
$(0, 1)$.  Each arc $K$ that is a component of
$f^{-1}(A)$ is mapped by $f$ in such a way that the endpoints of $K$
are mapped to the
same point of $C_2 \cap A$; call that point $a_j$.  We let $g$ agree
with $f$ on $\partial X$ and define $g$ on each component $K$ of
$f^{-1}(A)$ by defining $g(K) = a_j$ where $f(K) \subseteq A_j$ for
$j = 3, 4$.  Each component $V$ of $X - (f^{-1}(A) \cup \partial X)$
is a disk bounded by components of $f^{-1}(A)$ and arcs in $\partial X$.
Thus $g$ has been defined to send the simple closed curve $\partial
V$ to $C_2$ and it is
homotopic to $f$ on $\partial V$ where $f$ extends to $V$.  Since the
inclusion of $C_2$ into $X$ induces a monomorphism of the fundamental
groups, the degree of $g$ on $\partial V$ is zero and therefore it
extends to a map of $V$ to $C_2$.
In this way we define $g
\co X \to C_2$ as we claimed.
If $f(\partial C) = C_1$, we can apply the first part of the proof to
$\phi f \phi$ to obtain $g \co X \to C_2$.  By Lemma 6, the map
$\phi g \phi \co X \to C_1$ establishes that
$N_\partial(f) = MF_\partial[f]$ in this case.
\endproof

\section{Boundary semi-essential maps}

\begin{lemma} If $f \co (X, \partial X) \to  (X, \partial X)$
is a boundary semi-essential map where $f(C_1) = C_1$ and $f_1$ is
essential, then
$N_\partial(f) = MF_\partial[f]$.
\end{lemma}

\proof We assume $f$ is in reduced form then,
since $f_2$ is inessential, either
$f(C_2) = x_1 \in T_1$ or $f(C_2) = x_2 \in T_2$ so $f^{-1}(A)$ does
not intersect $C_2$.  By Theorem 6, there are
just four possible types of arcs that can be components of
$f^{-1}(A)$ and only those of type $[1, 3]$ fail to intersect $C_2$.
By Theorem 4, the only type of simple closed curve that could
be a component of $f^{-1}(A)$ is $(0, 1)$.  Referring
to Figure 6, the complement of the line segment $[y_1,
y_2]$ has two components that we call $U$ and $L$, where $C_2 \subset
L$.   All the simple closed curves in $f^{-1}(A)$ are contained in the
component of the complement of $h_2$ that contains $C_2$.   We have
chosen a point $z$ in $U \cap C_1$ such that the line segment $[x_0,
z]$ does not intersect $f^{-1}(A)$.  Since all the points of
$f^{-1}(A) \cap C_1$ lie in $U$, we may modify $f$ so that it
maps all of $L \cap C_1$ to $z$, without adding any fixed points to 
$f_1$.
The points $x_3$ and $x_4$ on the line segment between $x_0$ and $x_2$ 
are
chosen so that the segments $[x_0, x_3]$ and $[x_2, x_4]$ are
disjoint from $f^{-1}(A)$. We will construct a
map $g \co (X, \partial X) \to  (X, \partial X)$ that is homotopic
to $f$ as a map of pairs and identical to $f$ on the boundary of $X$
such that $g$ has no fixed points on the interior of $X$.   We
obtain $g$ on $L$ by defining it on the line segment $[x_0,
x_2]$ and then, using annulus extensions, we extend it to all of $L$.
We introduce notation as
follows.   Letting $h = h_1 \cup h_2$ we set $a_i = h \cap A_i$.  We
number the simple closed curves of $f^{-1}(A)$ as $K_1, \dots , K_m$
where $K_{j+1}$ is in the component of $X - K_j$ that contains $C_2$
and we set $k_j = K_j \cap [x_0, x_2]$.  We define $g \co [x_0,
x_2] \to [z, x_0] \cup [x_0, x_2] \cup h$ in the following way.  Let
$g(x_0) = z$ and extend to $[x_0, x_3]$ as a homeomorphism onto $[z,
x_0]$.  For
$j = 1, \dots, m$ we have $f(K_j) \subseteq A_{i(j)}$ for some $1 \le
i(j) \le 4$ and we define $g(k_j) = a_{i(j)}$.  We then define $g$ on
$[x_3, k_1]$ by sending it to the arc in $h$ connecting $x_0$ to
$a_{i(1)}$ and on each $[k_j, k_{j+1}]$ by mapping to an arc in $h$
that connects $a_{i(j)}$ to $a_{i(j+1)}$ through the component of $X -
A$ that contains $f(k_j, k_{j+1})$.   Setting $g(x_4) = x_0$, we
define $g$ on $[k_m, x_4]$ by mapping onto the arc in $h$ that
connects $a_{i(m)}$ to $x_0$.  Finally, $g$ maps $[x_4, x_2]$
homeomorphically onto $[x_0, f(C_2)]$.  We note that $g$ will have a 
fixed
point at $x_2$ in the case that $f(C_2) =
x_2$ and no fixed points if $f(C_2) = x_1$.  Choose two
simple closed curves in $L$ in the
complement of the $K_i$ intersecting $[x_0, x_2]$ in $x_3$ and $x_4$
respectively.  If we delete these simple closed curves as well as
all the $K_i$ from $L$, then the closure of each component is an
annulus containing a segment of $[x_0, x_2]$.
Applying annulus extensions to each of these, we 
now have $g$ defined on $L$.

\begin{figure}[ht!]
\centering
\includegraphics[scale=0.4]{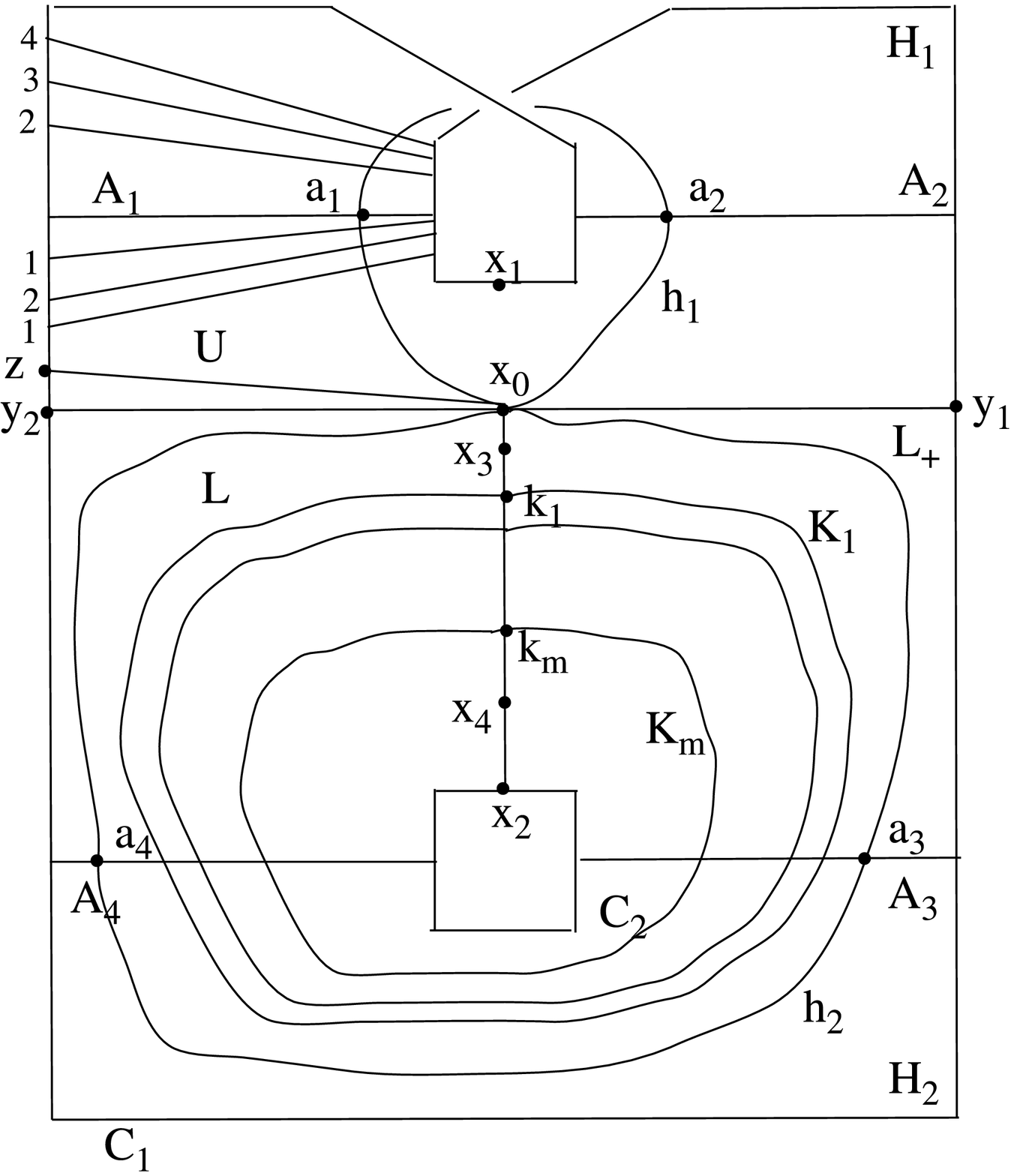}
\caption{} 
\end{figure}

The map $g$ is defined on $U$ as follows.  Each component of
$f^{-1}(A_3 \cup A_4)$ has both its endpoints mapped a single point,
either $C_1 \cap A_3$ or $C_1 \cap A_4$.  Consider a component $K$ of
$f^{-1}(A _1 \cup A_2)$ that is adjacent to a component $K'$ of
$f^{-1}(A_3 \cup A_4)$.  As the image of the two arcs of $C_1$
between $K$ and $K'$ must map to the same arc of $D \cap C_1$, then
the endpoints of $K$ must map to the same point, either $T_3 \cap
A_1$ or $T_4 \cap A_2$.  Next, for a component of $f^{-1}(A)$
adjacent to either $K$ or $K'$, its endpoints must also map to a
single point by the same argument.  Continuing in this manner, we
conclude that both endpoints of any component K of $f^{-1}(A)$ are
sent by $f$ to the same point $a_K \in A \cap C_1$.   Letting $g(K) =
a_K$ for each component $K$ of $f^{-1}(A)$ then, noting that $g([y_1,
y_2]) = z$, we can extend to all of $U$ so that $g(U) \subseteq C_1$
and $g^{-1}(A) = f^{-1}(A)$.  Since all the fixed points of $g$ in
$U$ are therefore in $C_1$, we conclude that $g$ has no fixed points
except those on $\partial X$ which are the fixed points of $f$ on
$\partial X$.  This completes the definition of $g$ on $X$ and, since
it has no fixed points in the interior of $X$, we have demonstrated
that $N_\partial(f) = MF_\partial[f]$.
\endproof

\begin{theorem} If $f \co (X, \partial X) \to  (X, \partial X)$
is a boundary semi-essential map such that $Im_\partial(f) = 2$, then
$N_\partial(f) = MF_\partial[f]$.
\end{theorem}

\proof The proof is divided into four cases, depending on
whether or not the boundary components are mapped to themselves and
which of the two boundary maps is essential.
The case in which each boundary
component is mapped
to itself by $f$, that $f_1$ is essential and $f_2$
is inessential is immediate from Lemma 8.
Next we consider the case where still $f_1$ is essential and $f_2$ is
inessential but now $f(C_1) = C_2$ and $f_2(C_2) = x_1$.
Since $f_2$ is inessential,
the components of $f^{-1}(A)$ are still the ones pictured in
Figure 6.   However, the arcs of
type $[1, 3]$ would now be labelled $3434\dots$ because only $A_3$ and 
$A_4$
intersect $C_2$.   Clearly,
both endpoints of each arc of type $[1, 3]$ must
map to the same point.  Thus we can define a map $g$ such that $g(U)
\subseteq C_2$ and $g([y_1, y_2]) = x_2$.  Our
definition of $g$ on $L$ will again consist of defining it on
$[x_0, x_2]$ and extending to all of $L$ by annulus extensions.
First, $g$ maps $[x_0,
x_3]$ homeomorphically onto $[x_2, x_0]$ with $g(x_3) = x_0$, then
$g$ will have exactly one fixed point on  $[x_0, x_3]$.  Next $g$
sends $[x_3, k_1]$ to the arc in $h$ that connects $x_0$ to $f(k_1)$.
As before, we continue in this manner so that $g$ maps $[k_1, x_4]$
to $h$.  Finally, $[x_4, x_2]$ is mapped homeomorphically onto the
line segment $[x_0, x_1]$.  Extending $g$ to all of $L$ by annulus
extensions, we note that the map $g$ has
only a single fixed point on all of $X$.  Since the Lefschetz number
of $f$ is nonzero, we conclude that $N_\partial(f) = MF_\partial[f] =
1$ in this case.

The remaining cases are (1) $f_1$ is an inessential map of $C_1$
to itself and $f_2$ is an essential map on $C_2$ and (2) $f_1$ is an
inessential map from $C_1$ to $C_2$ and $f_2$ is an essential map
from $C_2$ to $C_1$.  Both cases follow by Lemma 6 from
cases we have already established.
\endproof

The two results that follow complete the proof of Theorem 2 by
analyzing the cases in which $f$ is boundary semi-essential with
$Im_\partial(f) = 1$.

\begin{theorem} If $f \co (X, \partial X) \to  (X, \partial X)$
is a boundary semi-essential map such that $f(\partial X) = C_1$ and
$f_1$ is essential or if $f(\partial X) = C_2$ and $f_2$ is
essential, then
$N_\partial(f) = MF_\partial[f]$.
\end{theorem}

\proof   The case in which $f(\partial X) = C_1$ and $f_1$ is
essential is included in Lemma 8 and the other case then follows by 
Lemma 6.
\endproof

\begin{theorem} If $f \co (X, \partial X) \to  (X, \partial X)$
is a boundary semi-essential map such that $f(\partial X) = C_2$ and
$f_1$ is essential or if $f(\partial X) = C_1$ and $f_2$ is
essential, then
$ MF_\partial[f] - N_\partial(f) \le 1$.
\end{theorem}

\proof We will consider only the case that $f(\partial X) =
C_2$ and $f_1$ is essential since the other case will follow from
it by Lemma 6.   In contrast to the two previous
theorems, we cannot apply Lemma 8.  However, the fact that $f_1$
is essential still allows us to use Theorem 6 to reduce the possible
types of arcs among the components of $f^{-1}(A)$ to the four types
listed there.  Moreover, since $f_2$ is inessential, we may assume
that no point of $C_2$ is mapped to $A$ and therefore the only type
of arc is $[1, 3]$.  Theorem 4 thus assures us that we may still use
Figure 6 for the present case.  We will construct a map $g$ that is
boundary-preserving homotopic to $f$ with the same number of fixed
points on $\partial X$, that is $N(\bar f)$ of them, and one fixed
point in the interior of $X$.  The construction is a modification of
the proof of Lemma 8, as follows.  We may homotope $f$ on $C_1$
without changing the fixed points but mapping all of $C_1 \cap L$ to
$x_2$.  Define $g \co [x_0, x_2] \to [x_0, x_2] \cup h$ by
sending $[x_0, x_3]$ to $x_0$ and defining $g$ on $[x_3, x_2]$
exactly as in the proof of Lemma 8.  Extending $g$ to $L$ by annulus
extensions gives us a map without any
fixed points on the interior other than $x_0$.

In order to define $g$ on $U$, observe that both endpoints of a
component of $f^{-1}(A)$ in $U$ are mapped to the same point of $C_2
\cap A$.  Let $V$ be the component of $U - f^{-1}(A)$ containing
$[y_1, y_2]$, then $g$ can map $U - V$ to $C_2$ in the same way it
mapped $U$ to $C_1$ in the proof of Lemma 8.  The image under $g$ of
the boundary of $V$ is contained in $[x_0, x_2] \cup T_2$.
Since this set is contractible, we may extend $g$
from $\partial V$ to $V$ and thus $g$ is defined on all of $X$.
There are no fixed points of $g$ in $U$ except $x_0$ and those on
$C_1$.  We conclude that $g$ has a single fixed point on the interior 
of $X$.
\endproof

We conclude this paper with an observation about the difference
between\break $MF_\partial[f]$ and $N_\partial(f)$ for $f$ a
boundary-preserving map of the punctured M\"obius band $X$.   An
example in \cite{k2} shows that these are not always equal.  In this 
paper we
studied the question of equality by considering a classification of
boundary-preserving maps based on their restrictions to the two
boundary components of $X$.  From Theorems 7 - 12 we can see that the
only kind of map $f \co (X, \partial X) \to (X, \partial X)$ for
which $N_\partial(f)$ does not equal $MF_\partial[f]$ is one like the
example from \cite{k2}, which shows that Theorem 12 cannot be improved.

\Addresses
\end{document}